\documentclass{elsarticle}

\usepackage{amsfonts}
\usepackage{amsthm}
\usepackage{amsmath}
\usepackage{amssymb}
\usepackage{bbold}
\usepackage{mathtools} 
\usepackage{xcolor}
\usepackage{soul} 
\usepackage{mathrsfs}
\usepackage{tabularx}
\usepackage[top=1in, bottom=1in, left=1in, right=1in]{geometry}

\newtheorem{theorem}{Theorem}[section]
\newtheorem{lemma}{Lemma}[section]

\theoremstyle{definition}
\newtheorem{definition}{Definition}[section]

\theoremstyle{remark}

\newtheorem{example}{Example}[section]

\numberwithin{equation}{section}

\begin{document}
	
		\begin{frontmatter}
	
	\title{Second-order maximum principle controlled  weakly singular  Volterra integral equations}
	
	

    \author[]{Jasarat J. Gasimov}
    \ead{jasarat.gasimov@emu.edu.tr}
    \author[]{Nazim I. Mahmudov}
    \ead{nazim.mahmudov@emu.edu.tr}
    \cortext[cor1]{Corresponding author}

	\address{Department of Mathematics, Eastern Mediterranean University, Mersin 10, 99628, T.R. North Cyprus, Turkey}
	
	

	\begin{abstract}
	This work studies a class of singular Volterra integral equations that are (controlled) and can be applied to memory-related problems.For optimum controls, we prove a second-order Pontryagin type maximal principle. 
	\end{abstract}
\begin{keyword}
Second-order maximum principle, singular Volterra integral equation, optimal control, Pontryagin's maximum principle.
\end{keyword}  
	
		\end{frontmatter}
	\section{Introduction.}
Investigate the controlled Volterra integral equation given below:
	\begin{align}{\label{1}}
		y(t)=\eta(t)+\int_{0}^{t}\frac{f(t,s,y(s),u(s))}{(t-s)^{1-\alpha}}ds, \quad a.e.\quad t\in[0,T].
	\end{align}	
The above representations of $\eta(\cdot)$ and $f(\cdot,\cdot,\cdot,\cdot)$ are maps that are referred to as the generator and the free term of the state equation, respectively; $y(\cdot)$ represents the state trajectory and takes values in the Euclidean space $R^{n}$; $u(\cdot)$ represents the control and takes values in some separable metric space $U$. We present the cost functional performance metric to gauge the control's effectiveness.
\begin{align}{\label{2}}
	J(u(\cdot))=\int_{0}^{T}g(t,y(t),u(t))dt+\sum_{i=1}^{m}h^{i}(y(t_{i})),
\end{align}
The running cost and the prespecified instant costs $(at\quad 0\leq t_{1}<t_{2}<\cdots< t_{m}\leq T),$ are represented by the two terms on the right hand, respectively.

Memory-related dynamics can be described using equations such as \eqref{1}.Conversely, fractional (order) differential equations have garnered the interest of several academics in recent decades because of their intriguing applications in the fields of physics, chemistry, engineering, population dynamics, finance, and other different fields. Firstly, for the optimum control issue of singular integral equations, Pontryagin's type maximal principle was demonstrated by Lin, P and Yong, J in \cite{ly}, it was shown in \cite{mj} that there is a Pontryagin maximum principle for terminal state-constrained optimum control problems of Volterra integral equations with singular kernels. We want to make sure to point out that research was done on the maximum principle for stochastic equation optimum control issues \cite{ga}\cite{rbh},\cite{w},\cite{wy},\cite{wyy},\cite{hy},\cite{mb},\cite{an},\cite{yong},\cite{wang}. For the integral necessary condition of optimality of the second order for control problems given by volterra integral equations and a system of integro-differential equations, we refer to \cite{mk},\cite{ak},\cite{gam}, and for singular controls for systems with  fractional derivatives, and dynamic systems, see \cite{ym},\cite{yus},\cite{mt}

This research aims to investigate an optimum control problem of type \eqref{1} for a singular Volterra integral equation. We shall construct a second-order Pontryagin maximal principle with respect to the optimum control issue.
\section{Main Result.}
$\bullet$(A1): Let $f:\Delta \times R^{n}\times U\rightarrow R$ be a transformation with $(t,s)\mapsto f(t,s,y,,u)$ being measurable, $(y,u)\mapsto f(t,s,y,u)$ being continuously differentiable up to order-2. There is a constant $L>0$ such that 
\begin{align*}
	\vert f\vert+\vert f_{x}\vert+\vert f_{u}\vert\leq L(1+\vert x\vert+\vert u\vert).
\end{align*}

$\bullet$(A2): Let $g: [0,T]\times R^{n}\times U\rightarrow R$ be a transformation with $t\mapsto g(t,y,,u)$ being measurable, $(y,u)\mapsto g(t,y,u)$ being continuously differentiable up to order-2. There is a constant $L>0$ such that 
\begin{align*}
	\vert g\vert+\vert g_{x}\vert+\vert g_{u}\vert\leq L(1+\vert x\vert+\vert u\vert).
\end{align*}
\begin{definition}
	An admissible control $u(t)$ is considered to be singular according to the Pontryagin maximum principle if, within the process $\{u(t),y(t)\}$,  it implies that
		\begin{align}\label{H1}
		H_{u}(t) = 0, \quad t\in [0,T].
	\end{align}
\end{definition}
\begin{theorem}{\label{te}}
	Let $(A1)$ and $(A2)$ hold. Let $\eta(\cdot)\in L(0,T;R^{n})$ and $\eta(\cdot)$ be continuous at $t_{j}, j=1,2,\cdots,m.$ Suppose $(y^{*}(\cdot),u^{*}(\cdot))$ is an optimal pair of $\eqref{1}-\eqref{2}$.
	Then there a solution $\psi(\cdot)\in L^{\frac{p}{p-1}}(0,T;R^{n})$ of the following adjoint equation 
	\begin{align}{\label{4}}
		&\psi(t)=\int_{t}^{T}\frac{f_{y}(t,s,y^{*}(s),u^{*}(s))}{(s-t)^{1-\alpha}}\psi(s)ds-g_{y}(t,y^{*}(t),u^{*}(t))\\
		-&\sum_{i=1}^{m}1_{[0,t_{i})}(t)\frac{f_{y}(t_{i},t,y^{*}(t),u^{*}(t))}{(t_{i}-t)^{1-\alpha}}h^{i}_{y}(y^{*}(t_{i})), \quad a.e.\quad t\in[0,T]\nonumber,
	\end{align}
	such that following estimation holds:
\begin{align*}
\int_{0}^{T}H_{uu}(t)v^{2}(t)dt+\bigg[\int_{0}^{T}\int_{0}^{T}v(\tau)M(\tau,s)v(s)dsd\tau+2\int_{0}^{T}\bigg[\int_{0}^{t}v(s)H_{yu}(t)Q(t,s)v(t)dt\bigg]ds\bigg]\leq0.
\end{align*}
	where
	\begin{align}{\label{6}}
		&H(t,y,u,\psi)=\int_{t}^{T}\psi(s)\frac{f(t,s,y^{*}(s),u^{*}(s))}{(s-t)^{1-\alpha}}ds-g(t,y^{*}(t),u^{*}(t))\nonumber\\
		-&\sum_{i=1}^{m}1_{[0,t_{i})}(t)\frac{f(t_{i},t,y^{*}(t),u^{*}(t))}{(t_{i}-t)^{1-\alpha}}h^{i}_{y}(y^{*}(t_{i})).
	\end{align}
\end{theorem}
An extended Gronwall$^{,}$s inequality with a singular kernel is given in the following lemma.
\begin{lemma}\label{lem1}\cite{ly}
	Let $\alpha\in(0,1)$ and $q>\frac{1}{\alpha}$. Let $L(\cdot),a(\cdot),y(\cdot)$ be nonnegative functions with $L(\cdot)\in L^{q}(0,T)$ and $a(\cdot),y(\cdot)\in L^{\frac{q}{q-1}}(0,T)$. Suppose 
	\begin{align*}
		y(t)\leq a(t)+\int_{0}^{t}\frac{L(s)y(s)}{(t-s)^{1-\alpha}}ds, \quad a.e.\quad t\in [0,T].
	\end{align*}
Then there exists a constant $K>0$ such that 
\begin{align*}
	y(t)\leq a(t)+K\int_{0}^{t}\frac{L(s)a(s)}{(t-s)^{1-\alpha}}ds, \quad a.e.\quad t\in [0,T].
\end{align*}
\end{lemma}
\bigskip
Let $p\geq0$ and cosider  the following linear integral equation:
\begin{align}{\label{fr}}
y(t)=\eta(t)+\int_{0}^{t}\frac{A(t,s)y(s)}{(t-s)^{1-\alpha}}ds, \quad a.e.\quad t\in[0,T].
\end{align}	
	where $\alpha\in(0,1), \eta(\cdot)\in L^{p}(0,T;R^{n})$, and $A:\Delta\to R^{n\times n}$ is measurable and satisfies
	\begin{align*}
		\vert A(t,s)\vert\leq L(s), \quad (t,s)\in \Delta,
	\end{align*}
for some measurable function $L(\cdot)\in L^{(\frac{1}{\alpha}\vee \frac{p}{p-1})_{+}}(0,T)$.
\begin{lemma}\label{lem2}\cite{ly}
	Let $1\leq p<\frac{1}{1-\alpha}$, for any  $s\in[0,T) $
	\begin{align*}
		\Phi(t,s)=\frac{A(t,s)}{(t-s)^{1-\alpha}}+\int_{s}^{t}\frac{A(t,\tau)\Phi(\tau,s)}{(t-\tau)^{1-\alpha}}d\tau,\quad a.e.\quad t\in(s,T]. 
	\end{align*}
Then
\begin{align*}
	y(t)=\eta(t)+\int_{0}^{t}\Phi(t,s)\eta(s)ds,\quad a.e.\quad \in[0,T]
\end{align*}
the expression gives a representation for the solution to the linear equation \eqref{fr}.
\end{lemma}
\subsection{Proof of the theorem.}
Let $(y^{*}(\cdot),u^{*}(\cdot))$ be an optimal pair of $\eqref{1}-\eqref{2}$ and fix any $u(\cdot)\in U^{p}[0,T].\\$
Denote
\begin{align}{\label{ig}}
	u^{\delta}(\cdot)=u{\cdot}+\delta v(\cdot)\quad where\quad v(\cdot)=u(\cdot)-u^{*}(\cdot).
\end{align}
Clearly, $u^{\delta}\in U^{p}[0,T]$. Let $y^{\delta}(\cdot)=y(\cdot,\eta(\cdot),u^{\delta}(\cdot))$ be the corresponding state.

It follows that
\begin{align}{\label{7}}
	&y^{\delta}(t)-y^{*}(t)=\int_{0}^{t}\frac{f(t,s,y^{\delta}(s),u^{\delta}(s))-f(t,s,y^{*}(s),u^{*}(s))}{(t-s)^{1-\alpha}}ds\\
	=&\int_{0}^{t}\frac{f^{\delta}_{y}(t,s)}{(t-s)^{1-\alpha}}(y^{\delta}(t)-y^{*}(t))ds+\int_{0}^{t}\frac{f^{\delta}_{u}(t,s)}{(t-s)^{1-\alpha}}(u^{\delta}(t)-u^{*}(t))ds\nonumber
\end{align}
where 
\begin{align}{\label{15}}
	&f^{\delta}_{y}(t,s)=\int_{0}^{1}f_{y}(t,s,y^{*}(s)+\tau[y^{\delta}(t)-y^{*}(t)],u^{\delta}(s))d\tau, \quad (t,s)\in\Delta,\nonumber\\
	&f^{\delta}_{u}(t,s)=\int_{0}^{1}f_{u}(t,s,y^{\delta}(s),u^{*}(s)+\tau[u^{\delta}(t)-u^{*}(t)])d\tau, \quad (t,s)\in\Delta.
\end{align}
(A1) provides with
\begin{align}{\label{17}}
	\vert f^{\delta}_{y}(t,s)\vert\leq L(s), \quad \vert f^{\delta}_{u}(t,s)\vert\leq L(s), \quad (t,s)\in\Delta.
\end{align}
Clearly, $L(s)\in L^{q}(0,T)$ for some $q\in (\frac{1}{\alpha},p)$. That being so
\begin{align}
	&\vert y^{\delta}(t)-y^{*}(t)\vert\nonumber\\
	=&\int_{0}^{t}\frac{L(s)}{(t-s)^{1-\alpha}}\vert y^{\delta}(t)-y^{*}(t)\vert ds+\int_{0}^{t}\frac{L(s)}{(t-s)^{1-\alpha}}\vert u^{\delta}(t)-u^{*}(t)\vert ds, \quad t\in[0,T].
\end{align}
By the extended Gronwall$^{,}$s inequality Lemma \ref{lem1} and \eqref{ig}, choosing $q^{\prime}\in(\frac{1}{\alpha},q)$(see, \cite{ly}),
\begin{align}\label{42}
	\vert y^{\delta}(t)-y^{*}(t)\vert\leq K\delta^{\frac{q-q^{\prime}}{q^{\prime}q}}\to0, \delta\to0 \quad uniformly\quad in\quad t\in[0,T].
\end{align}

$\bullet$ let $Y_{1}(\cdot)$ is the solution of the following first-order variational equation:
\begin{align}{\label{8}}
	&Y_{1}(t)=\int_{0}^{t}\frac{f_{y}(t,s,y^{*}(s),u^{*}(s))}{(t-s)^{1-\alpha}}Y_{1}(s)ds+\int_{0}^{t}\frac{f_{u}(t,s,y^{*}(s),u^{*}(s))}{(t-s)^{1-\alpha}}v(s)ds\quad t\in[0,T].
\end{align}
$\bullet$ let $Y_{2}(\cdot)$ is the solution of the following second-order variational equation:
\begin{align}{\label{100}}
	&Y_{2}(t)=\int_{0}^{t}\frac{f_{y}(t,s,y^{*}(s),u^{*}(s))}{(t-s)^{1-\alpha}}Y_{2}(s)ds+\int_{0}^{t}\frac{f_{yy}(t,s,y^{*}(s),u^{*}(s))}{(t-s)^{1-\alpha}}Y^{2}_{1}(s)ds\nonumber\\
	+&\int_{0}^{t}\frac{2f_{yu}(t,s,y^{*}(s),u^{*}(s))}{(t-s)^{1-\alpha}}Y_{1}(s)v(s)ds+\int_{0}^{t}\frac{f_{uu}(t,s,y^{*}(s),u^{*}(s))}{(t-s)^{1-\alpha}}v^{2}(s)ds t\in[0,T].
\end{align}
As a consequence of those,
\begin{align}{\label{9}}
&\frac{y^{\delta}(t)-y^{*}(t)}{\delta}-Y_{1}(t)\nonumber\\
=&\int_{0}^{t}\frac{f^{\delta}_{y}(t,s)}{(t-s)^{1-\alpha}}\bigg(\frac{y^{\delta}(t)-y^{*}(t)}{\delta}\bigg)ds+\int_{0}^{t}\frac{f^{\delta}_{u}(t,s)}{(t-s)^{1-\alpha}}\bigg(\frac{u^{\delta}(t)-u^{*}(t)}{\delta}\bigg)ds\nonumber\\
=&\int_{0}^{t}\frac{f_{y}(t,s,y^{*}(s),u^{*}(s))}{(t-s)^{1-\alpha}}Y_{1}(s)ds+\int_{0}^{t}\frac{f_{u}(t,s,y^{*}(s),u^{*}(s))}{(t-s)^{1-\alpha}}v(s)ds\\
=&\int_{0}^{t}\frac{f^{\delta}_{y}(t,s)}{(t-s)^{1-\alpha}}\bigg(\frac{y^{\delta}(t)-y^{*}(t)}{\delta}-Y_{1}(s)\bigg)ds+\int_{0}^{t}\frac{f^{\delta}_{u}(t,s)}{(t-s)^{1-\alpha}}\bigg(\frac{u^{\delta}(t)-u^{*}(t)}{\delta}-v(s)\bigg)ds\nonumber\\
+&\int_{0}^{t}\frac{f^{\delta}_{y}(t,s)-f_{y}(t,s,y^{*}(s),u^{*}(s))}{(t-s)^{1-\alpha}}Y_{1}(s)ds+\int_{0}^{t}\frac{f^{\delta}_{u}(t,s)-f_{u}(t,s,y^{*}(s),u^{*}(s))}{(t-s)^{1-\alpha}}v(s)ds\nonumber\\
=&\int_{0}^{t}\frac{f^{\delta}_{y}(t,s)}{(t-s)^{1-\alpha}}\bigg(\frac{y^{\delta}(t)-y^{*}(t)}{\delta}-Y_{1}(s)\bigg)ds+a^{\delta}_{1}(t)+a^{\delta}_{2}(t)+a^{\delta}_{3}(t)\quad  t\in[0,T].\nonumber
\end{align}
The dominated convergence theorem (like \cite{ly}), and \eqref{ig} supplies
\begin{align*}
	\lim_{\delta\to0}\vert a^{\delta}_{1}(t)\vert=0,\quad
	\lim_{\delta\to0}\vert a^{\delta}_{2}(t)\vert=0,\quad
	\lim_{\delta\to0}\vert a^{\delta}_{3}(t)\vert=0.
\end{align*}
The dominated convergence theorem and the extended Gronwall$^{,}$s inequality Lemma \ref{lem1}(like \cite{ly}) produces
\begin{align}\label{41}
		\lim_{\delta\to0}\bigg\vert \frac{y^{\delta}(t)-y^{*}(t)}{\delta}-Y_{1}(t)\bigg\vert=0.
\end{align}
Sequentially,
\begin{align*}
	&\frac{y^{\delta}(t)-y^{*}(t)}{\delta}-Y_{1}(t)-\frac{\delta}{2}Y_{2}(t)\\
	=&\int_{0}^{t}\frac{f_{y}(t,s,y^{*}(s),u^{*}(s))}{(t-s)^{1-\alpha}}\bigg(\frac{y^{\delta}(s)-y^{*}(s)}{\delta}-Y_{1}(s)-\frac{\delta}{2}Y_{2}(s)\bigg)ds\nonumber\\
	+&\frac{f_{u}(t,s,y^{*}(s),u^{*}(s))}{(t-s)^{1-\alpha}}\bigg(\frac{u^{\delta}(s)-u^{*}(s)}{\delta}-v(s)\bigg)ds\\
	+&\frac{1}{2\delta}\int_{0}^{t}\frac{f_{yy}(t,s,y^{*}(s),u^{*}(s))}{(t-s)^{1-\alpha}}(y^{\delta}(s)-y^{*}(s))^{2}ds-\frac{\delta}{2}\int_{0}^{t}\frac{f_{yy}(t,s,y^{*}(s),u^{*}(s))}{(t-s)^{1-\alpha}}Y^{2}_{1}(s)ds\\
	+&\frac{1}{\delta}\int_{0}^{t}\frac{f_{yu}(t,s,y^{*}(s),u^{*}(s))}{(t-s)^{1-\alpha}}(y^{\delta}(s)-y^{*}(s))(u^{\delta}(s)-u^{*}(s))ds-\frac{\delta}{2}\int_{0}^{t}\frac{f_{yu}(t,s,y^{*}(s),u^{*}(s))}{(t-s)^{1-\alpha}}Y_{1}(s)v(s)ds\\
	+&\frac{1}{2\delta}\int_{0}^{t}\frac{f_{uu}(t,s,y^{*}(s),u^{*}(s))}{(t-s)^{1-\alpha}}(u^{\delta}(s)-u^{*}(s))^{2}ds-\frac{\delta}{2}\int_{0}^{t}\frac{f_{uu}(t,s,y^{*}(s),u^{*}(s))}{(t-s)^{1-\alpha}}v^{2}(s)ds\quad  t\in[0,T].\nonumber
\end{align*}
The extended Gronwall$^{,}$s inequality Lemma \ref{lem1}, \eqref{41}, and \eqref{ig} yields
\begin{align}\label{40}
	&\lim_{\delta\to0}\bigg\vert \frac{y^{\delta}(t)-y^{*}(t)}{\delta}-Y_{1}(t)-\frac{\delta}{2}Y_{2}(t)\bigg\vert=0.
\end{align}
Also, by the optimality of $(y^{*}(\cdot),u^{*}(\cdot))$, one has
\begin{align}{\label{10}}
&0\leq J(u^{\delta}(\cdot))-J(u^{*}(\cdot))=\int_{0}^{T}[g(t,y^{\delta}(t),u^{\delta}(t))-g(t,y^{*}(t),u^{*}(t))]dt\nonumber\\
=&\int_{0}^{T}g_{y}(t,y^{*}(t),u^{*}(t))(y^{\delta}(t)-y^{*}(t))dt+\int_{0}^{T}g_{u}(t,y^{*}(t),u^{*}(t))(u^{\delta}(t)-u^{*}(t))dt\\
+&\frac{1}{2}\int_{0}^{T}g_{yy}(t,y^{*}(t),u^{*}(t))(y^{\delta}(t)-y^{*}(t))^{2}dt+\int_{0}^{T}g_{yu}(t,y^{*}(t),u^{*}(t))(y^{\delta}(t)-y^{*}(t))(u^{\delta}(t)-u^{*}(t))dt\nonumber\\
+&\frac{1}{2}\int_{0}^{T}g_{uu}(t,y^{*}(t),u^{*}(t))(u^{\delta}(t)-u^{*}(t))^{2}dt+\sum_{i=1}^{m}h^{i}_{y}(y^{*}(t_{i}))(y^{\delta}(t)-y^{*}(t))+\frac{1}{2}\sum_{i=1}^{m}h^{i}_{yy}(y^{*}(t_{i}))(y^{\delta}(t)-y^{*}(t))^{2}.\nonumber
\end{align}
Further, from \eqref{40},\eqref{41}, and \eqref{42}, we get
\begin{align}{\label{11}}
	&0\leq J(u^{\delta}(\cdot))-J(u^{*}(\cdot))=\delta\int_{0}^{T}g_{y}(t,y^{*}(t),u^{*}(t))Y_{1}(t)dt+\frac{\delta^{2}}{2}\int_{0}^{T}g_{y}(t,y^{*}(t),u^{*}(t))Y_{2}(t)dt\nonumber\\
	+&\delta\int_{0}^{T}g_{v}(t,y^{*}(t),u^{*}(t))v(t)dt+\frac{\delta^{2}}{2}\int_{0}^{T}g_{yy}(t,y^{*}(t),u^{*}(t))Y^{2}_{1}(t)dt\\
	+&\delta^{2}\int_{0}^{T}g_{yu}(t,y^{*}(t),u^{*}(t))Y_{1}(t)v(t)dt+\frac{\delta^{2}}{2}\int_{0}^{T}g_{uu}(t,y^{*}(t),u^{*}(t))v(t)^{2}dt\nonumber\\
	+&\sum_{i=1}^{m}h^{i}_{y}(y^{*}(t_{i}))Y_{1}(t_{i})+\frac{\delta^{2}}{2}\sum_{i=1}^{m}h^{i}_{y}(y^{*}(t_{i}))Y_{2}(t_{i})+\frac{\delta^{2}}{2}\sum_{i=1}^{m}h^{i}_{yy}(y^{*}(t_{i}))Y^{2}_{1}(t_{i})+o(\delta^{2}),\quad (\delta\to0).\nonumber
\end{align}
Applying \eqref{4},\eqref{8},\eqref{100},\eqref{6}, and Fubini$^{\prime}$s theorem,
\begin{align}{\label{12}}
	&0\leq J(u^{\delta}(\cdot))-J(u^{*}(\cdot))=-\delta\int_{0}^{T}H_{u}(t)v(t)dt\nonumber\\
	-&\frac{\delta^{2}}{2}\int_{0}^{T}H_{yy}(t)Y^{2}_{1}(t)dt-\delta^{2}\int_{0}^{T}H_{yu}(t)Y_{1}(t)v(t)dt-\frac{\delta^{2}}{2}\int_{0}^{T}H_{uu}(t)v^{2}(t)dt\\
	+&\frac{\delta^{2}}{2}\sum_{i=1}^{m}h^{i}_{yy}(y^{*}(t_{i}))Y^{2}_{1}(t_{i})+o(\delta^{2}),\quad (\delta\to0).\nonumber
\end{align}
Lemma \ref{lem2} provides following expression
\begin{align}{\label{Ex}}
	Y_{1}(t)=\int_{0}^{t}Q(t,s)v(s)ds
\end{align}
Substitute \eqref{Ex} into \eqref{12}, and Fubini$^{\prime}$s theorem
\begin{align*}
	&\int_{0}^{T}H_{yy}(t)Y^{2}_{1}(t)dt=\int_{0}^{T}H_{yy}(t)\bigg(\int_{0}^{t}Q(t,s)v(s)ds\bigg)\bigg(\int_{0}^{t}Q(t,\tau)v(\tau)d\tau\bigg)dt\\
	=&\int_{0}^{T}\int_{0}^{T}v(s)\bigg[\int_{\max\{\tau,s\}}^{T}Q(t,s)H_{yy}(t)Q(t,\tau)dt\bigg]v(\tau)dsd\tau=\int_{0}^{T}\int_{0}^{T}v(s)M(\tau,s)v(\tau)dsd\tau,
\end{align*}
and
\begin{align*}
\sum_{i=1}^{m}h^{i}_{yy}(y^{*}(t_{i}))Y^{2}_{1}(t_{i})=\int_{0}^{T}\int_{0}^{T}v(s)\bigg[\sum_{i=1}^{m}1_{[0,t_{i})}(s)1_{[0,t_{i})}(\tau)Q(t_{i},s)h^{i}_{yy}(y^{*}(t_{i}))Q(t_{i},\tau)\bigg]v(\tau)d\tau ds.
\end{align*}
Therefore,
\begin{align}{\label{13}}
	&0\leq J(u^{\delta}(\cdot))-J(u^{*}(\cdot))=-\delta\int_{0}^{T}H_{u}(t)v(t)dt-\frac{\delta^{2}}{2}\int_{0}^{T}H_{uu}(t)v^{2}(t)dt\nonumber\\
	-&\frac{\delta^{2}}{2}\bigg[\int_{0}^{T}\int_{0}^{T}v(\tau)M(\tau,s)v(s)dsd\tau+2\int_{0}^{T}\bigg[\int_{0}^{t}v(s)H_{yu}(t)Q(t,s)v(t)dt\bigg]ds\bigg]+o(\delta^{2}),\quad (\delta\to0)
\end{align}
where
\begin{align*}
M(\tau,s)=\int_{\max\{\tau,s\}}^{T}Q(t,s)H_{yy}(t)Q(t,\tau)dt-\sum_{i=1}^{m}1_{[0,t_{i})}(s)1_{[0,t_{i})}(\tau)Q(t_{i},s)h^{i}_{yy}(y^{*}(t_{i}))Q(t_{i},\tau),\quad \tau,s\in[0,T].
\end{align*}
We can now prove our theorem, which is based on the cost functional estimation mentioned previously.
	
The theorem is obtained by dividing the right side of \eqref{13} by $\delta^{2}$ and allowing $\delta$ to approach zero, while taking into account $H_{u}(t) = 0 $ in expression (\ref{13}).

The theorem was proved.
\begin{example}
	Consider  the problem
	\begin{align*}
		y(t)=1+t\sqrt{t}+\int_{0}^{t}\frac{ty(s)u(s)}{(t-s)^{\frac{1}{2}}},\quad a.e \quad t\in[0,1],
	\end{align*}
\begin{align*}
	J(u)=y(1)+\int_{0}^{1}y(s)u(s)ds\longrightarrow min,\quad \vert u\vert\leq1.
\end{align*}
We are evaluating the efficiency of the control input $u(t) = 1$ and analyzing its optimality. This particular selection of control corresponds to the solution $1 + t\sqrt{t}$ for the integral equation. Throughout the course of the process represented by $(1 + t\sqrt{t}, 0)$, we have noted the following outcome.
\begin{align*}
	\psi(t)=1,\quad  H=0.
\end{align*}
Therefore, the control $u(t) = 0$ is identified as a singular control. Clearly, employing the control $u(t) = 0$ yields a performance measure value of $J(u) = 2$. Now, let's investigate if there is an alternative control function that leads to functional values less than $2$. We will compute the value of $J$ for the admissible control $u(t) = -\frac{1}{2}$.
\begin{align*}
	y(t)=1
\end{align*} 
Then, we have 
\begin{align*}
	\frac{1}{2}=J(-\frac{1}{2})\leq J(0)=2
\end{align*}
This indicates that opting for the control $u(t) = 0$ within the interval $t \in [0, 1]$ is not an optimal.
\end{example}

\end{document}